\documentclass[12pt]{amsart}
\usepackage{amsmath,amstext,amsfonts,amsthm,amssymb}
\usepackage{eucal}
\usepackage{diagrams}                                     %
\diagramstyle[scriptlabels,PostScript=dvips]               %
\newarrow{Dash}{}{dash}{}{dash}{}                          %
\renewcommand{\subsubsection}[1]{\addtocounter{subsubsection}{1}
{\ \\[3pt]\bf \thesubsubsection. \  #1} }

\swapnumbers

\newtheorem{CRL}[subsection]{Corollary}
\newtheorem*{Thm}{Theorem}

{  \theoremstyle{definition}

}

\newcommand{\cT}{\mathcal{T}}

\renewcommand{\P}{\mathbb{P}}

\newcommand{\C}{\mathbb{C}}

\newcommand{\fM}{\mathfrak{M}}
\newcommand{\fX}{\mathfrak{X}}
\newcommand{\eps}{\varepsilon}
\renewcommand{\phi}{\varphi}

\begin{document}
\title[Plumbing coordinates]
{Plumbing coordinates on Teichm\"uller space: a counterexample}
\author{Vladimir Hinich}
\address{Department of Mathematics, University of Haifa, Mount Carmel,
 Haifa 31905,  Israel}
\email{hinich@math.haifa.ac.il}

\begin{abstract}We present an example showing that a family of Riemann
surfaces obtained by a general plumbing construction
does not necessarily give local coordinates on the Teichm\"uller space.
\end{abstract}
\maketitle

\section{Introduction}
The following well-known construction gives a holomorphic family
of punctured Riemann surfaces of genus $g$ with $n$ punctures.

Fix a maximally degenerated Riemann surface $X_0$ of genus $g$ with $n$
punctures. The irreducible components of $X_0$ are triply-punctured Riemann
spheres $S_i$, $i=1,...,d=2g-2+n$; some pairs of punctures are glued
together to form $m=3g-3+n$ double points; and $n$ remaining punctures are
the punctures of $X_0$.

Choose local coordinates around $2m$ punctures corresponding to the
double points of $X_0$ by fixing $2m$ open embeddings $z_i,w_i:D\rTo X_0$,
$i=1,\ldots,m$,
of the unit disc $D=\{z||z|<1\}$ so that $z_i(0)$ and $w_i(0)$ correspond
to the same $i$-th double point of $X_0$. We assume the following
disjointness property for the embeddings
$z_i,\ w_j$.
\begin{itemize}
\item[(DIS)] {\em The images of the
open embeddings $z_i,\ w_j$ have no intersection, except for the double
points 
 $z_i(0)=w_i(0)$.}
\end{itemize}

The choice of the coordinates gives rise to a family $\fX=\{X_t\}$
of (nodal) genus $g$ Riemann surfaces with $n$ punctures, based on $D^m$;
for $t=(t_1,\ldots,t_m)$ the Riemann surface $X_t$ is obtained by
replacing the neighborhood of the $i$-th doble point given by the equation
$z_iw_i=0$ with the surface given by the equation $z_iw_i=t_i$.

The described above family can be presented as a flat morphism
$\pi:\fX\rTo D^m$ of analytic spaces, so that $X_t$ identifies with the
fiber $\pi^{-1}(t)$. The fiber $X_t$ is non-singular
for $t=(t_1,\ldots,t_m)$ iff all $t_i$ are nonzero.
Let $D_0=D-\{0\}$ and $B=D_0^m$. \newpage
We have therefore
a holomorphic map $\alpha$ from $B$ to the moduli stack $\fM_{g,n}$.
\footnote{
A choice of marking on the maximally degenerate curve $X_0$
gives a lifting of $\alpha$ to a map $\beta:B\rTo\cT_{g,n}/\Gamma$,
$\Gamma$ being the subgroup of the modular group generated by the
Dehn twists around the preimages of the double points. There is no
difference for us which one of the maps $\alpha$ or $\beta$ to consider.
}

In this note we show that the map $\alpha:B\rTo\fM_{g,n}$ is not necessarily
\'etale(=a local isomorphism). This seems to contradict Corollary on p. 346
of A.~Marden's account of his work with C.~Earle (see~\cite{marden}):

{\em For an arbitrary choice of local coordinates at the double points
of $X_0$, the map $B\rTo\cT_{g,n}/\Gamma$, $\Gamma$ being the subgroup
generated by the Dehn twists around the preimages of the double points,
 is a local isomorphism.
}



For the horocyclic coordinates at the double points, the claim of the above
statement follows from the theorem at p. 345 of Marden's paper~\cite{marden}
(see Kra~\cite{kra} for a detailed study of horocyclic coordinates).

As it was kindly pointed to us by the referee, this easily implies
the following weaker version of the claim which could have been
what Marden had really in mind:

\begin{CRL}\label{rect}
For an arbitrary choice of local coordinates at the double points of $X_0$,
the map $\beta:B\rTo\cT_{g,n}/\Gamma$ provides
local coordinates for the quotient space $\cT_{g,n}/\Gamma$
\emph{in a neighborhood of zero}.
\end{CRL}

The reasoning is as follows. Let $\beta_{EM}:B\rTo\cT_{g,n}/\Gamma$
be the map based on the horocyclic coordinates ($B$ is now a neighborhood
of $0$ in the punctured polydisc). The composition
$$ \beta_{EM}^{-1}\circ\beta:B\rTo B$$
extends to a holomorphic map of the polydiscs preserving the coordinate
hyperplanes. This implies that the Jacobian matrix of the map
at zero is diagonal. Since one circuit around $j$-th coordinate
corresponds in both polydiscs to the same Dehn twist, the diagonal entries
of the Jacobian matrix are nonzero. This implies Corollary~\ref{rect}.

Another proof of Corollary~\ref{rect}, based on the deformation
theory, can be found in~\cite{HV}, 5.3.7.

We say that a family of curves $\pi:\fX\rTo B$ provides a local chart
near $t\in B$ if the corresponding map $\alpha:B\rTo\fM_{g,n}$ is \'etale
at $t$.

Earle-Marden coordinates provide a local chart at any point (and even global
coordinates). However, these coordinates do not meet the disjointness
assumptions (DIS). This is sometimes inconvenient.

The following result of~\cite{HV}, Section 5.3.5--5.3.6, shows
that the the local coordinates based on the plumbing construction
satisfying (DIS) cover the whole Teichm\"uller space.
\begin{Thm}
Let $X\in\fM_{g,n}$.  Then there exists a choice of local coordinates
for disjoint neighborhoods of the punctures of the triply-punctured spheres,
 for which $X$ lies in the image of the local chart.
\end{Thm}

\section{The counterexample}
\subsection{Kodaira-Spencer calculation}
Let $\pi:\fX\rTo B$ be the family of Riemann surfaces constructed
from a choice of a maximally degenerate Riemann surface $X_0$ and a
collection of local coordinates around the punctures. Let $t\in B$
and put $X=X_t$.

We wish to describe the map of the tangent spaces induced by $\alpha$ at
$t$. The tangent space to $\alpha(t)$ is, according to Kodaira-Spencer,
the first cohomology $H^1(X,T)$ of $X$ with the coefficients
in the sheaf of holomorphic vector fields on $X$ vanishing at the punctures.

Let $\tau=(\tau_1,\ldots,\tau_m)$ be the tangent vector at $t\in B$.
It is easy to write down the \v{C}ech cocycle corresponding to $\tau$.
The curve $X$ is covered by the pairs of pants $P_i$ obtained from the
triply-punctured sphere $S_i$ by cutting discs $\{z||z|<|t_j|\}$ if $j$-th
double point corresponds to a puncture of $S_i$. Each double point in $X_0$
corresponds to an annulus $A_j$ in $X_t$ defined in the coordinates $z_j$
or $w_j$ by the inequality $|t_j|<|z_j|<1$ or, what is the same,
 $|t_j|<|w_j|<1$. Using the formulas (4.10) of \cite{Kod}, 4.2, we get
that the tangent vector $\tau$ is sent to the cohomology class defined
by the \v{C}ech cocycle assigning to each annulus $A_j$ the vector
field
\begin{equation}
\label{cocycle}
\frac{\tau_j}{t_j}z_j \frac{d}{dz_j}.
\end{equation}

It is not, however, very easy to decide, in general, whether the cocycle
(\ref{cocycle}) defines a non-zero element in the cohomology.

\subsection{The case $g=0,\ n=4$}

We  will construct a counterexample in the first nontrivial case
$g=0, \ n=4$. Here $m=1$ and $X_0$ is a union of two triply-punctured
spheres. In this case $X$ is a sphere with four punctures and there is
a very easy way of determining non-triviality of a one-cocycle.
In fact, for $X=\P^1_\C\setminus\{a,b,c,\infty\}$ the only regular quadratic
differential has form
\begin{equation}\label{eq:q}
q=dw^2/(w-a)(w-b)(w-c)
,
\end{equation}
where $w$ is the standard coordinate in $\P^1_\C$.
Then a cocycle is non-trivial iff its pairing with $q$ is nonzero.

In order to make explicit calculations, we will change the setup.

Let $L=\{z|1-\eps<|z|<1+\eps\}$. Let $\phi:L\rTo\P^1_\C$ be a holomorphic
open embedding, so that $\P^1_\C\setminus\phi(L)$ consists of two components,
$X$ and $Y$.
In what follows we use the variabe $w$ to denote the standard coordinate
in  $\P^1_\C$.
 Let $a,b,c\in\C$ such that  that $a,b\in X$ and $c,\infty\in Y$.
Put $\widetilde{X}=X\cup\phi(L)$ and $\widetilde{Y}=Y\cup\phi(L)$.

The Riemann sphere $(\P^1_\C;a,b,c,\infty)$ with four punctures can be obtained
by a plumbing construction from two triply punctured spheres
as follows. The first sphere $\P_1$ is obtained by gluing
the disc $\{z|1+\eps>|z|\}$ to $\widetilde{X}$ via the identification
$z\sim\phi(z)$ for $z\in L$. The punctures of $\P_1$ are $w=a,b$ in
$\widetilde{X}$ and the point $z=0$.
The second sphere $\P_2$ is obtained by gluing the disc
$\{z|1-\eps<|z|\leq\infty\}$ to $\widetilde{Y}$ via the same identification
$z\sim\phi(z),\ z\in L.$ The punctures of $\P_2$ are $w=c,\infty$ in
$\widetilde{Y}$ and the point $z=\infty$.
The plumbing coordinates are $z_1=\frac{z}{1+\epsilon}$ for $\P_1$
and $w_1=\frac{1-\epsilon}{z}$ for $\P_2$. For the value of plumbing parameter
$t=\frac{1-\epsilon}{1+\epsilon}$ the annuli $t<|z_1|<1$ and $t<|w_1|<1$
identify with $L$ and give $\P^1_\C$ with the punctures $w=a,b,c,\infty$ as
required.

Note that our presentation of $(\P^1_\C,a,b,c,\infty)$ as a result of plumbing
construction does not allow to get explicit formulas for the plumbing
parametra $z_1,w_1$ in terms of the standard coordinates in $\P_1$ and $\P_2$.
But this is the price we pay for having an explicit coordinate $w$ in the
resulting four-punctured sphere. The intersection of two pairs of pants
covering our Riemann surface is the image $\phi(L)$ of the annulus $L$.
The vector field on $\phi(L)$ is proportional to $z\frac{d}{dz}$ in the
$z$-coordinate of $L$. Thus, its pairing with the nontrivial quadratic
differential~(\ref{eq:q}) is given by the integral
\begin{equation}
\label{integral}
\int_C\frac{z(\phi'(z))^2dz}{(\phi(z)-a)(\phi(z)-b)(\phi(z)-c)},
\end{equation}
the integral being taken along the unit circle.

\subsection{The counterexample}
We wish to find $\phi,a,b,c$ such that the integral~(\ref{integral}) vanishes.
We will be looking for $\phi$ given by the formula
\begin{equation}
\label{formula-for-phi}
\phi(z)=(z-R)^2,
\end{equation}
where $R\in\C$ is outside the unit circle, so that the function $\phi$,
restricted to the unit circle, is an open embedding.

A way to ensure the relative position of $a,b,c$ with respect to $\phi(L)$
is the following. We will choose $a=\phi(A), \ b=\phi(B), \ c=\phi(C)$,
so that
\begin{equation}
\label{conditions}
|A|<1,\ |B|<1, A\not=B,\ |C|>1,\ |R|>1, \ |C-2R|>1.
\end{equation}
The integral~(\ref{integral}) can be now easily calculated. The result is
\begin{equation}
\int_C\frac{4z(z-R)^2dz}{(z-A)(z+A-2R)(z-B)(z+B-2R)(z-C)(z+C-2R)}
\end{equation}
which is equal to, up to the constant $4\pi i$,
\begin{equation}
\frac{1}{(A-B)(A+B-2R)}\left(\frac{A(A-R)}{(A-C)(A+C-2R)}-
 \frac{B(B-R)}{(B-C)(B+C-2R)}\right).
\end{equation}
Thus, the integral vanishes if
\begin{equation}
A(A-R)(B-C)(B+C-2R)=B(B-R)(A-C)(A+C-2R).
\end{equation}
Opening the brackets and dividing by $A-B$, we get
\begin{equation}
A=\frac{C(R-B)(C-2R)}{C^2-2RC+RB}.
\end{equation}
We will show now the equation has a solution satisfying the constraints
(\ref{conditions}).
Choose $B,R$ such that $|B|<1,\ |R|>1$ but $|R-B|<1$. If we make $C$
tending to infinity, $A$ with tend to $R-B$.

This proves there exist values of the parameters satisfying the constraints
and making the integral vanish.

\subsection{Acknowledgement} The author is very grateful to Zalman
Rubinstein for a good advice. The referee's suggestions have helped very
much to improve the exposition.

\end{document}